\documentstyle[amssymb,12pt]{article}

\def\nabc{\displaystyle{\buildrel c\over \nabla}}
\def\oisk{\buildrel\textstyle{\rm k \atop ^\vee} \over\dots,}
\def\oisr{\buildrel\textstyle{\rm r \atop ^\vee} \over\dots}
\def\oiss{\buildrel\textstyle{\rm s \atop ^\vee} \over\dots}
\def\Talpha#1{\vbox{\ialign{##\crcr
    $\alpha$\crcr\noalign{\kern2pt\nointerlineskip}
           $\hfil\displaystyle{#1}\hfil$\crcr}}} \def\Out{\mbox{Out}}
\def\Outp{\mbox{\scriptsize Out}} \def\Int{\mbox{Int}}

\def\Onabla#1{\vbox{\ialign{##\crcr
$\,\scriptstyle{0}$\crcr\noalign{\kern2pt\nointerlineskip}
           $\hfil\displaystyle{#1}\hfil$\crcr}}} 
  
 \def\oplusinf{\mathop{\oplus}}
\def\otimesinf{\mathop{\otimes}}  \def\calc{{\cal C}}
 \def\k{\Bbb K} \def\bbbone{\mbox{\rm 1\hspace
{-.6em} l}}  \def\hom{{\mbox{Hom}}}
 \def\der{{\mbox{\scriptsize Der}}}
\def\gder{{\mbox{Der}}} 
\def\mer{{\mbox{\normalshape\scriptsize Der}}}
\def\os{\underline{\Omega}} \def\hs{\underline{H}}
\def\gmer{{\mbox{\normalshape Der}}} 

 \newtheorem{theorem}{THEOREM}
\newtheorem{lemma} {LEMMA}

\begin{document}

\baselineskip=0.7cm \begin{center} {\Large\bf CONNECTIONS ON CENTRAL
BIMODULES} \end{center} \vspace{0.75cm}

\begin{center} Michel DUBOIS-VIOLETTE \\ \vspace{0.3cm} {\small
Laboratoire de Physique Th\'eorique et Hautes
Energies\footnote{Laboratoire associ\'e au Centre National de la
Recherche Scientifique - URA D0063}\\ Universit\'e Paris XI, B\^atiment
211\\ 91 405 Orsay Cedex, France\\ E-mail: flad@qcd.th.u-psud.fr}\\
\vspace{0.2cm} and\\ \vspace{0.2cm} Peter W. MICHOR\\ \vspace{0.3cm}
{\small Erwin Schr\"odinger Institute of Mathematical Physics\\
Pasteurgasse 6/7\\ A-1090 Wien, Austria\\ E-mail: michor@esi.ac.at}
\end{center} \vspace{1cm}

\begin{center} March 29, 1995 \end{center}

\vspace {1cm}

\noindent L.P.T.H.E.-ORSAY 94/100\\ \noindent ESI-preprint 210\\
\noindent q-alg/9503020

\newpage \begin{abstract} We define and study the theory of
derivation-based connections on a recently introduced class of bimodules
over an algebra which reduces to the category of modules whenever the
algebra is commutative. This theory contains, in particular, a
noncommutative generalization of linear connections. We also discuss the
different noncommutative versions of differential forms based on
derivations. Then we investigate reality conditions and a noncommutative
generalization of pseudo-riemannian structures.  \end{abstract}

\section{Introduction and notations}

There are several noncommutative generalizations of the calculus of
differential forms and, more generally, of the differential calculus of
classical differential geometry, e.g. [2 to 10]. As stressed in [3], the
extension of classical tools to the noncommutative setting is never
straightforward. This means that, in order to produce relevant objects,
one must have in mind a lot of examples coming both from mathematics and
from physics. In this paper, we concentrate on the differential calculus
based on derivations as generalization of vector fields, [4]. It was
shown in [5] that this differential calculus is natural for quantum
mechanics in the sense that with it, quantum mechanics has the same
relation to noncommutative symplectic geometry as classical mechanics to
classical symplectic geometry.  For finite quantum spin systems this was
already pointed out in [6].\\ In this paper, $A$ is an associative
algebra over $\k=\Bbb R$ or $\Bbb C$ with a unit $\bbbone$. The algebra
$A$ is to be considered as the generalization of the algebra of smooth
functions and the Lie algebra $\gder(A)$ of all derivations of $A$ as
the generalization of the Lie algebra of smooth vector fields. The Lie
algebra $\gder(A)$ is also a module over the center $Z(A)$ of $A$ and
furthermore $Z(A)$ is stable by the action of $\gder(A)$. The
corresponding Lie algebra homomorphism of $\gder(A)$ into the Lie
algebra $\gder(Z(A))$ factorizes through the Lie algebra $\Out(A)$ of
all derivations of $A$ modulo the ideal $\Int(A)$ of all inner
derivations of $A$; the Lie algebra $\Out(A)$ is also a $Z(A)$-module.
Notice that if $A$ is commutative, $A=Z(A)$ and $\gder(A)=\Out(A)$; so
$\Out(A)$ is also a generalization of the Lie algebra of vector fields
and this is a good generalization for a theory of ``invariants". Indeed
in general one has $H^0(A,A)=Z(A)$ and $H^1(A,A)=\Out(A)$, (whereas
$\gder(A)=Z^1(A,A))$, where $H(A,A)$ is the Hochschild cohomology of $A$
with value in $A$. So $Z(A)$ and $\Out(A)$ are Morita invariant as well
as the homomorphism of $\Out(A)$ into $\gder(Z(A))$.  We now recall the
relevant generalizations of differential forms in this context [4], [9].
As for the commutative case [11], the notions of differential forms can
be extracted from the differential algebra $C(\gder(A),A)$ of
Chevalley-Eilenberg cochains of the Lie algebra $\gder(A)$ with values
in the $\gder(A)$-module $A$. There are two natural generalizations of
the graded differential algebra of differential forms which use
$\gder(A)$ as generalization of vector fields : A minimal one,
$\Omega_{\der}(A)$, which is the smallest differential subalgebra of
$C(\gder(A),A)$ which contains $A$ and a maximal one, $\os_{\der}(A)$,
which consists of all cochains in $C(\gder(A),A)$ which are
$Z(A)$-multilinear.\\ As mentioned above, it is also useful to use
$\Out(A)$ as generalization of vector fields. The corresponding
generalizations of differential forms $\Omega_{\Outp}(A)$ and
$\os_{\Outp}(A)$ are respectively graded differential subalgebras of
$\Omega_{\der}(A)$ and $\os_{\der}(A)$. To obtain them, one notices that
there is a canonical operation, in the sense of H. Cartan [1], $X\mapsto
i_X$ for $X\in \gder(A)$, of the Lie algebra $\gder(A)$ in the graded
differential algebra $C(\gder(A),A)$ defined by
$i_X\alpha(X_1,\dots,X_{n-1})=\alpha(X,X_1,\dots,X_{n-1})$ for $X_k\in
\gder(A)$ and $\alpha\in C^n(\gder(A),A)$. Both $\Omega_{\der}(A)$ and
$\os_{\der}(A)$ are stable by the $i_X,X\in \gder(A)$, and
$\Omega_{\Outp}(A)$ and $\os_{\Outp}(A)$ are defined to be the
respective differential subalgebras which are basic with respect to the
corresponding operation of $\Int(A)$, i.e. one has:
$$\Omega_{\Outp}(A)=\{\alpha\in \Omega_{\der}(A)\vert i_X\alpha=0\
\mbox{and}\ L_X\alpha =0,\ \forall X\in \Int(A)\}$$
$$\os_{\Outp}(A)=\{\alpha\in \os_{\der}(A)\vert i_X\alpha=0\ \mbox{and}\
L_X\alpha=0,\ \forall X\in \Int(A)\},$$ where $L_X=di_X + i_Xd$ as
usual. One has the inclusions of graded differential algebras \[
\begin{array}{lll} \Omega_{\der}(A)& \subset & \os_{\der}(A)\\ \bigcup &
& \bigcup\\ \Omega_{\Outp}(A) & \subset & \os_{\Outp}(A) \end{array} \]
In the case where $A$ is the algebra of smooth functions on a
finite-dimensional paracompact smooth manifold, all these graded
differential algebras coincide with the graded differential algebra of
differential forms. In general, there is a differential calculus for $A$
in $\Omega_{\der}(A)$ and in $\os_{\der}(A)$. However if $A$ is not
commutative, i.e. $A\not=Z(A)$, then $\Omega_{\Outp}(A)$ and
$\os_{\Outp}(A)$ do not contain $A$ and are not $A$-modules. So they do
not carry a differential calculus for $A$. The differential algebra
$\os_{\Outp}(A)$ can be identified with the differential algebra
$C_{Z(A)}(\Out(A),Z(A))$ of $Z(A)$-linear cochains of the Lie algebra
$\Out(A)$ with values in $Z(A)$. So $\os_{\Outp}(A)$ is a Morita
invariant generalization of differential forms. We shall use
$\os_{\der}(A)$ for the differential calculus and then, the
``invariants" will be closed elements in the subalgebra $\os_{\Outp}(A)$
leading to Morita-invariants in the cohomology $\hs_{\Outp}(A)$.\\ In
this paper, we wish to extend, for $A$ noncommutative, the theory of
connections (derivation laws) on $A$-modules for $A$ commutative as
formulated in [11]. There are several noncommutative generalizations of
the notion of module over a commutative algebra. First one can consider
the notion of right (or left) $A$-module. Alternatively, one can
remember that a module over a commutative algebra is canonically a
bimodule of a very specific kind and we speak of the {\sl induced
structure of bimodule}. In [8], we introduced the notion of {\sl central
bimodule}: This is just a $A$-bimodule such that the underlying
structure of $Z(A)$-bimodule is induced by a structure of $Z(A)$-module,
i.e. multiplication by elements of $Z(A)$ on both sides coincide. This
notion is stable by arbitrary projective and inductive limits and by
tensor products over $A$ or over $Z(A)$. When $A$ is commutative, a
central bimodule is just a module (for the induced bimodule structure).
It is for this notion that we define and study connections in this
paper. There are several reasons to prefer this notion rather than that
of right or left module. The first one is that our one-forms constitute
such a bimodule and that we wish to be able to define linear
connections. A second very general reason, which is connected with
quantum mechanics, is explained in the remark of Section 8. In [8] and
[9] we also introduced the more restrictive notion of {\sl diagonal
bimodule}: This is a bimodule isomorphic to a subbimodule of $A^I$, for
some set $I$, where $A$ is equipped with its canonical structure of
$A$-bimodule. A diagonal bimodule is central and, if $A$ is commutative,
a diagonal bimodule is just a module such that the canonical mapping
into its bidual is injective. Both $\Omega_{\der}(A)$ and
$\os_{\der}(A)$ are diagonal and therefore central; this is why the
notion of connection considered here includes a generalization of the
notion of linear connection. Furthermore, and this was the very reason
diagonal bimodules were introduced, it was shown in [8] that the
derivation (differential) $d:A\rightarrow \Omega^1_{\der}(A)$ is
universal for derivations of $A$ into diagonal bimodules: i.e. for any
derivation $\delta:A\rightarrow M$ of $A$ into a diagonal bimodule $M$,
there is a unique bimodule homomophism
$i_\delta:\Omega^1_{\der}(A)\rightarrow M$ such that
$\delta=i_\delta\circ d$.\\ Finally we shall need, to describe torsion
for instance, the generalization of vector valued differential forms. It
was shown in [9] that the right spaces to generalize the
Fr\"olicher-Nijenhuis bracket were the space $\gder(A,\Omega_{\der}(A))$
of derivations of $A$ into $\Omega_{\der}(A)$ if one uses
$\Omega_{\der}(A)$ as generalization of differential forms and the space
$\gder(A,\os_{\der}(A))$ if one uses $\os_{\der}(A)$ as generalization
of differential forms. In this paper it is this latter generalization
that will be considered. If $N$ and $M$ are $A$-bimodules, we use the
notation $\hom^A_A(N,M)$ to denote the space of bimodule homomorphisms
of $N$ into $M$. This is a $Z(A)$-bimodule which is in fact a
$Z(A)$-module whenever $M$ is central.\\ The plan of the paper is the
following. In Section 2 we define the notion of derivation-based
connection on central bimodules. In Section 3 we describe some
constructions which allow to produce new connections from given
connections. In Section 4 we define linear connections and their
torsions. In Section 5 we give some basic examples. In Section 6 we
introduce and study a duality between bimodules and modules over the
center. In Section 7 we apply this duality to the one-forms showing, in
particular, that $\os_{\der}^1(A)$ is the bidual of $\Omega_{\der}^1(A)$
for this duality. In Section 8 we study reality conditions for the case
of $\ast$-algebras. Finally, in Section 9 we investigate a
noncommutative generalization of pseudo-riemannian structures in our
framework.

\section{Connections on central bimodules}

Let $M$ be a central bimodule over $A$, {\sl a connection on} $M$ is a
linear mapping $\nabla$, $X\mapsto \nabla_X$, of $\gder(A)$ into the
linear endomorphisms of $M$ such that one has \[ \left\{
\begin{array}{l} \nabla_{zX}(m)=z\nabla_X(m)\\ \nabla_X(amb)=X(a)mb +a
\nabla_X(m) b + am X(b) \end{array} \right.  \] $\forall m\in M, \forall
X \in \gder(A), \forall z\in Z(A)$ and $\forall a,b\in A$.\\ Given
$\nabla$ as above, {\sl the curvature} $R$ {\sl of} $\nabla$ is the
bilinear antisymmetric mapping $(X,Y)\mapsto R_{X,Y}$ of $\gder(A)
\times \gder(A)$ into the linear endomorphisms of $M$ defined by
$$R_{X,Y}(m)=\nabla_X(\nabla_Y(m))-\nabla_Y(\nabla_X(m))-\nabla_{[X,Y]}(m),$$
$$\forall X,Y\in \gder A,\ \forall m\in M.$$ One has the following
properties \[ \left\{ \begin{array}{l} R_{zX,Y}(m)=zR_{X,Y}(m),\\
R_{X,Y}(amb)=aR_{X,Y}(m) b \end{array} \right.\] $\forall m\in M,
\forall X,Y\in \gder(A), \forall z\in Z(A), \forall a,b\in A.$\\ Thus
$R$ is an antisymmetric $Z(A)$-bilinear mapping of $\gder(A)\times
\gder(A)$ into the $Z(A)$-module $\hom^A_A(M,M)$ i.e. $$R\in
\hom_{Z(A)}\ (\Lambda^2_{Z(A)}\gder(A), \hom^A_A(M,M)).$$ From its very
definition and from the Jacobi identity, it follows that $R$ satisfies
the {\sl Bianchi identity} \[
[\nabla_X,R_{Y,Z}]+[\nabla_Y,R_{Z,X}]+[\nabla_Z,R_{X,Y}]= R_{[X,Y],Z} +
R_{[Y,Z],X}+R_{[Z,X],Y}.  \] There is another way to describe all that.
Let $\os^n_{\der}(A,M)$ be the space (in fact the $Z(A)$-module) of
antisymmetric $Z(A)$-multilinear mappings of $(\gder(A))^n$ into $M$,
i.e. one has $$\os^n_{\der}(A,M)=\hom_{Z(A)}(\Lambda^n_{Z(A)}\gder
(A),M).$$ The spaces $\os^n_{\der}(A,M)$ as well as
$$\os_{\der}(A,M)=\oplusinf_n \os^n_{\der}(A,M)$$ are canonically
A-bimodules which are central bimodules. Then a connection $\nabla$ as
above on $M$ is simply a linear mapping of $M$ into $\os^1_{\der}(A,M)$
which satisfies $$\nabla(amb)= \displaystyle{da\otimesinf_A}
mb+a\nabla(m)b+\displaystyle{am\otimesinf_A} db,\ \forall a,b \in A\
\mbox{and}\ \forall m\in M,$$ where the canonical bimodule homomorphisms
$$\os^1_{\der}(A)\displaystyle{\otimesinf_A} M \rightarrow
\os^1_{\der}(A,M)\ \mbox{and}\ M\displaystyle{\otimesinf_A}
\os^1_{\der}(A) \rightarrow \os^1_{\der}(A,M)$$ have been used.\\ More
generally, by using the canonical bimodule homomorphisms
$$\os^m_{\der}(A)\displaystyle{\otimesinf_A}
\os^n_{\der}(A,M)\rightarrow \os^{m+n}_{\der}(A,M)$$ and
$$\os^n_{\der}(A,M)\displaystyle{\otimesinf_A}
\os^m_{\der}(A)\rightarrow \os^{m+n}_{\der}(A,M),$$ one equips
$\os_{\der}(A,M)$ with a structure of graded $\os_{\der}(A)$-bimodule.
Let us extend $\nabla:\os^0_{\der}(A,M)\rightarrow \os^1_{\der}(A,M)$ to
an endomorphism, again denoted by $\nabla$, of $\os_{\der}(A,M)$ with
$\nabla (\os^n_{\der}(A,M))\subset \os^{n+1}_{\der}(A,M)$ by the
following definition

\[ \begin{array}{ll} (\nabla \varphi)(X_0,\dots,X_n)&
=\displaystyle{\sum_{0\leq k \leq n}}(-1)^k
\nabla_{X_k}\varphi(X_0,\oisk,X_n)\\ & +\displaystyle{\sum_{0\leq
r<s\leq n}}(-1)^{r+s}\varphi([X_r,X_s],X_0,\oisr\oiss,X_n) \end{array}
\]
 for $\varphi\in\os^n_{\der}(A,M)$ and $X_k \in \gder (A)$, where
$\buildrel \textstyle{\rm k \atop ^\vee}\over .$ means omission of
$X_k$. One has, for $\alpha\in \os^a_{\der}(A)$,
$\beta\in\os^b_{\der}(A)$ and $\varphi\in\os^n_{\der}(A,M)$:
$$\nabla(\alpha\varphi\beta)=(d\alpha)\varphi\beta+(-1)^a\alpha\nabla(\varphi)\beta+(-1)^{a+n}\alpha\varphi
d\beta.$$ It follows that $\nabla^2$ which is the canonical extension of
the curvature satisfies
$\nabla^2(\alpha\varphi\beta)=\alpha\nabla^2(\varphi)\beta$, i.e. it is
a homomorphism of $\os_{\der}(A)$-bimodules (and of graded
$\os_{\der}(A)$-bimodules) of $\os_{\der}(A,M)$ into itself, (the
Bianchi identity now reads $\nabla(\nabla^2)=(\nabla^2)\nabla)$.

\section{Associated connections}

There exist central bimodules which do not admit connections. For
instance, in [11], J.L. Koszul gives the following example: take $A=\Bbb
K[t]$, i.e. the commutative algebra of polynomials in $t$, and $M=A/N$
where $N$ is the ideal of polynomials without constant term; then $M$ is
a central bimodule since it is an $A$-module with $A$ commutative and
there is no connection on $M$ because if $\nabla$ is such a connection
and if $e$ denotes the class of $\bbbone$ in $A/N$, one must have
$$0=\nabla_{\partial/\partial t}(te)=e+t\nabla_{\partial/\partial
t}(e)=e,$$ i.e. a contradiction. However, if $X\mapsto \nabla_X$ is a
connection on a central bimodule $M$ and if $X\mapsto \Gamma_X$ is a
$Z(A)$-linear mapping of $\gder(A)$ into $\hom^A_A(M,M)$ then $X\mapsto
\nabla_X + \Gamma_X$ is also a connection on $M$ and all connections on
$M$ are of this form; i.e. if the set of connections on a central
bimodule $M$ is not empty, it is an affine space modelled on
$\hom_{Z(A)}(\gder(A),\hom^A_A(M,M))$. Notice that, for $M=A$,
$\nabla_X(a)=X(a)\ (\forall a\in A,\ \forall X\in \gder(A))$ is a
connection on $A$ with vanishing curvature which will be referred to as
the {\sl canonical connection on} $A$. In this section, we will describe
connections on central bimodules associated with bimodules which admit
connections. These connections will be accordingly called {\sl
associated connections}.\\ Let $M$ be a central bimodule equipped with a
connection $\nabla$ and let $N$ be a subbimodule of $M$. Assume that
$\nabla_XN\subset N$ for any $X\in \gder(A)$. Then the restriction of
$\nabla$ to $N$, ( i.e. of the $\nabla_X,\ X\in \gder(A))$, is a
connection on $N$ and $\nabla$ induces a connection on the quotient
bimodule $M/N$. In both cases, we shall speak of the {\sl induced
connections} by $\nabla$ to design these connections on $N$ and on
$M/N$.\\ Let $(M_i)_{i\in I}$ be a family of central bimodules equipped
with connections $\nabla^i$. Then $\nabla_X((m_i)_{i\in
I})=(\nabla^i_X(m_i))_{i\in I}$, for $m_i\in M_i$ and $X\in\gder(A)$,
defines a connection on the product $\displaystyle{\prod_{i\in I}}M_i$.
By restriction, one obtains a connection on the direct sum
$\displaystyle{\oplusinf_{i\in I}}M_i$,
$\nabla_X(\sum_im_i)=\sum_i\nabla^i_X m_i$, since
$\nabla_X(\displaystyle{\oplusinf_{i\in I}}M)\subset
\displaystyle{\oplusinf_{i\in I}}M_i$ for any $X\in \gder(A)$. These
connections will be called {\sl product} and {\sl direct sum} of the
connections $\nabla^i$. One defines similarily {\sl projective limits}
and {\sl inductive limits} of connections when the appropriate stability
conditions are satisfied.\\ Let $M$ and $M'$ be two central bimodules
equipped with connections $\nabla$ and $\nabla'$. For $X\in \gder(A)$,
consider the linear endomorphisms $\nabla_X\otimes id_{M'}+id_M\otimes
\nabla'_X$ of $M\otimes M'$. The bimodule $M\otimes M'$ is not central
in general, however the subbimodules generated, respectively by the $$ma
\otimes m'-m\otimes am',\ a\in A,\ m\in M,\ m'\in M'$$ and by the
$$mz\otimes m'-m\otimes zm',\ z\in Z(A),\ m\in M,\ m'\in M$$ are stable
by the above endomorphisms (remembering that\linebreak[4]
$\gder(A)(Z(A))\subset Z(A))$, so they define endomorphisms of
$M\displaystyle{\otimesinf_A}M'$ and of
$M\displaystyle{\otimesinf_{Z(A)}}M'$ which are easily seen to be
 connections on $M\displaystyle{\otimesinf_A}M'$ and
$M\displaystyle{\otimesinf_{Z(A)}}M'$, respectively. These connections
will be called
 {\sl tensor product} of $\nabla$ and $\nabla'$ over $A$ and over
$Z(A)$, respectively. By induction, one defines the tensor product (over
$A$ or over $Z(A))$ of a finite family of connections on a finite family
of central bimodules. This tensor product is associative in an obvious
sense.\\ In particular, if $M$ is a central bimodule with a connection
$\nabla$, then by applying the above construction, one obtains a
connection $\nabla^\otimes$ on the tensor algebra of $M$ over $A$,
$T_A(M)=\displaystyle{\oplus_n}(\otimes^n_A M)$, satisfying
$\nabla^\otimes_X(a)=X(a)$ for $a\in A=T^0_A(M)$ and $X\in \gder(A)$.
One has $\nabla^\otimes_X(tt')=\nabla^\otimes_X(t) t' +
t\nabla^\otimes_X(t')$ for $t,t'\in T_A(M),\ X\in \gder(A)$.\\ Let $M$
be a central bimodule, then $\hom^A_A(M,M)$ is an algebra over $Z(A)$.
The group of invertible elements of $\hom^A_A(M,M)$ will be called {\sl
the group of gauge transformations of} $M$. Given a connection $X\mapsto
\nabla_X$ and a gauge transformation $g$ on $M$, $X\mapsto g
\circ\nabla_X\circ g^{-1}$ is again a connection which will be referred
to as {\sl the gauge transform of} $\nabla$ {\sl by} $g$. Two
connections belonging to the same orbit will be called {\sl gauge
equivalent connections}.

\section{The case $M=\os^1_{\mbox{$\der$}}(A)$: Linear connections}

The bimodule $\os^1_{\der}(A)$ is diagonal and therefore central. A
connection $\nabla$ on $\os^1_{\der}(A)$ will be called a {\sl linear
connection}. There is a canonical bimodule homomorphism
$\mu:\os^1_{\der}(A,\os^1_{\der}(A))\rightarrow \os^2_{\der}(A)$ which
extends the product $\os^1_{\der}(A)\displaystyle{ \otimesinf_A}
\os^1_{\der}(A)\rightarrow \os^2_{\der}(A)$, namely
$\mu(\varphi)(X,Y)=\varphi_X(Y)-\varphi_Y(X)$ for $X,Y\in \gder(A)$ and
$\varphi\in \os^1_{\der}(A,\os^1_{\der}(A))$. Given a linear connection
$\nabla$, one defines a linear mapping $T$ of $A$ into $\os^2_{\der}(A)$
by setting $T(a)=-\mu \circ \nabla (da)$ for $a\in A$. One has
$T(ab)=T(a)b + aT(b)$ for $a,b\in A$, therefore $T$ is an element of
$\gder(A,\os^2_{\der}(A))$ which will be called {\sl the torsion} of the
linear connection $\nabla$. Since $\os^2_{\der}(A)$ is a diagonal
bimodule, it follows from the universal property of the derivation
$d:A\rightarrow \Omega^1_{\der}$ that there is a unique bimodule
homomorphism $i_T:\Omega^1_{\der}(A)\rightarrow \os^2_{\der}(A)$ such
that $T=i_T\circ d$. The explicit form of $i_T$ is easy to write, one
has $i_T=d-\mu\circ \nabla$ which extends as a bimodule homomorphism of
$\os^1_{\der}(A)$ into $\os^2_{\der}(A)$. We shall frequently identify
the torsion $T\in \gder(A,\os^2_{\der}(A))$ with this element of
$\hom^A_A(\os^1_{\der}(A),\os^2_{\der}(A))$.

\section{Examples}

\subsection{The case where $A$ is commutative}

In the case where $A$ is commutative, a central bimodule is simply an
$A$-module and the notion of connection defined here reduces to the
usual one, i.e. to the notion of derivation laws of [11]. One obtains
the classical notion of connection on a smooth vector bundle $E$ of
finite rank over a smooth finite-dimensional paracompact manifold $V$ by
taking the algebra $C^\infty(V)$ of smooth functions on $V$ for $A$ and
by taking the module $\Gamma(E)$ of smooth sections of $E$, i.e. a
typical finite projective module over $A=C^\infty(V)$. Since the
canonical mapping of $\Gamma(E)$ into its bidual is injective, the
underlying bimodule is not only central but it is also a diagonal
bimodule.\\ Now we investigate cases which are of ``opposite side".

\subsection{The case where $\Out(A)=~0$}

Let us now assume that $A$ is a noncommutative algebra which has only
inner derivations, i.e. $\Int(A)=\gder(A)$ or, equivalently $\Out(A)=0$.
In this case, every central bimodule $M$ admits a canonical connection
$\nabc$ with vanishing curvature defined by: $\nabc_{ad(x)}(m)=xm-mx$,
$\forall x\in A$ and $\forall m \in M$. The other connections on $M$ are
of course of the form $\nabla_{ad(x)}=\nabc_{ad(x)}+\Gamma_{ad(x)}$
where $\Gamma\in \hom_{Z(A)}(\Int(A), \hom^A_A(M,M))$. Since the
curvature of $\nabc$ vanishes one cannot have a non trivial theory of
characteristic classes using the above notion of connection for such
algebras. This also partly explains why, in the general case, one has to
factorize the inner derivations out in order to get a good theory of
invariants.\\ For $M=\underline{\Omega}^1_{\der}(A)$, $\nabc$ is a
linear connection.  Its torsion $T$ is given by
$T(a)(ad(x),ad(y))=-ad[x,y](a)=-[[x,y],a]$, or
$i_T(\omega)(ad(x),ad(y))=-\omega(ad([x,y]))$, for $x,y,a\in A$,
$\omega\in \underline{\Omega}^1_{\der}(A)$.

\subsection{The case where $A$ has a trivial center $Z(A)=\k.\bbbone$}

In this case, $Z(A)$-linearity reduces to $\k$-linearity, so in
particular the Lie derivative $X\mapsto L_X=i_Xd+di_X$ is a connection
on any of the central bimodules $\Omega^n_\der(A)$ and
$\underline{\Omega}^n_\der(A)$. These connections have vanishing
curvatures since the Lie derivative is a homomorphism of Lie algebras.
Acting on $\underline{\Omega}^1_\der(A)$ the Lie derivative is then a
linear connection with a torsion $T$ given by $T(a)(X,Y)=[X,Y](a)$, or
$i_T(\omega)(X,Y)=\omega([X,Y])$, for $a\in A$, $X,Y\in \gder (A)$,
$\omega\in \underline{\Omega}^1_\der(A)$.\\ Notice that if one has also
$\Out(A)=0$, then both $\nabc$ and $L$ are connections with zero
curvature on the $\Omega^n_\der(A)$ and $\underline{\Omega}^n_\der(A)$
but in general they are not gauge equivalent, except for $n=0$ where
they coincide with the canonical connection on $A$. In particular, on
$\underline{\Omega}^1_\der(A)$ they are linear connections with opposite
torsion and therefore $\frac{1}{2}(\nabc+L)$ is (on
$\underline{\Omega}^1_\der(A))$ torsion-free.

\subsection*{Remarks} A priori, examples 5.2 and 5.3 are independent
(Morita invariant) classes of algebras. For instance if $C$ is a unital
commutative algebra which is different from ${\k}.\bbbone$ and which has
no nonzero derivation, e.g. $C={\k}^n$ with $n\geq 2$, then the matrix
algebra $M_N(C)$ has a non-trivial center, $C$, and all its derivations
are inner; on the other hand, if $E$ is a vector space of dimension
$\geq 2$, the tensor algebra $T(E)$ of $E$ has a trivial center but any
non vanishing endomorphism of $E$ extends uniquely as a derivation of
$T(E)$ which is never inner. However since here $A$ is the analog of the
algebra of smooth functions, one could prefer to choose $A$ in such a
way that it has ``many" derivations. From this point of view, it is
natural to introduce the following class $\calc^{\infty,0}$: $A$ belongs
to the class $\calc^{\infty,o}$ if $X(a)=0$, $\forall X\in \gder(A)$ for
$a\in A$ implies $a\in {\k}.\bbbone$. It is worth noticing here that
this condition might not be sufficient to ensure the existence of
``many" derivations: For instance let $A=\oplus A^n$ be a $\Bbb
Z$-graded algebra with $A^0= \k.\bbbone$, then the degree derivation
defined by $deg(a)=na$ if $a\in A^n$ is such that $deg(a)=0$ implies
$a\in {\k.\bbbone}$, so $A$ is in $\calc^{\infty,0}$ but it is easy to
construct examples such that the only derivations are the multiple of
$deg$. In any case, any $A$ in $\calc^{\infty,0}$ such that $\Out(A)=0$
has a trivial center (i.e.  examples 5.2 in $\calc^{\infty,0}$ are
contained in examples 5.3).

\section{Duality and diagonal bimodules}

Let $M$ be a central bimodule over $A$, then the space $\hom^A_A(M,A)$
of all bimodule homomorphisms of $M$ into $A$ is a module over the
center $Z(A)$ of $A$, i.e. it is a $Z(A)$-module which will be denoted
by $M^{\ast_{A}}$ and called the {\sl dual of the bimodule} $M$ when no
confusion arises. The reader must be aware of the fact that
$M^{\ast_{A}}$ {\sl is not} the dual of $M$ as $A\otimes A^{op}$-module
or as $A\displaystyle{\otimesinf_{Z(A)}}A^{op}$-module. Conversely, let
$N$ be a $Z(A)$-module then the space $\hom_{Z(A)}(N,A)$ is canonically
a bimodule over $A$ which is diagonal, and therefore central, since it
is a subbimodule of $A^N$. This diagonal bimodule will be denoted by
$N^{\ast_{A}}$ and called {\sl the dual of the} $Z(A)$-module $N$. Thus
one has a duality between central bimodules over $A$ and modules over
$Z(A)$ which obviously refers to $A$; this duality is similar to the
duality between left and right $A$-modules. In fact, when $A$ is
commutative {\sl all these four notions coincide} with the notion of
$A$-module. Notice that if $M$ is a central bimodule, {\sl the duality}
$(M,M^{\ast_{A}})$ {\sl is separated if and only if} $M$ {\sl is
diagonal}; another way to say the same thing is to remark that there is
a canonical bimodule homomorphism $c_M:M\rightarrow
M^{\ast_{A}\ast_{A}}$ and that this canonical homomorphism is injective
if and only if $M$ is diagonal. Dually, if $N$ is a $Z(A)$-module, then
there is a canonical $Z(A)$-module homomorphism $c_N:N\rightarrow
N^{\ast_{A}\ast_{A}}$ which is in general not injective nor surjective;
a sufficient condition for the injectivity of $c_N$ is that the
canonical mapping of $N$ into its $Z(A)$-module bidual
$N^{\ast_{Z(A)}\ast_{Z(A)}}$ is injective. A $Z(A)$-module $N$ will be
said to be $A$-{\sl diagonal}, or simply {\sl diagonal} if no confusion
arises, whenever the canonical mapping $c_N$ is injective or, which is
the same, whenever it is separated by $N^{\ast_{A}}=\hom_{Z(A)}(N,A)$;
this means that it is isomorphic to a $Z(A)$-submodule of $A^I$ for some
set $I$. Thus the dual $M^{\ast_{A}}$ of any central bimodule $M$ is
diagonal. More generally, {\sl a duality} between a central bimodule $M$
and a $Z(A)$-module $N$ will be a bimodule homomorphism
$\langle,\rangle$ of $M\displaystyle{\otimesinf_{Z(A)}}N$ into $A$,
$(m,n)\mapsto \langle m,n\rangle$; the duality $\langle,\rangle$ {\sl is
separated in} $M$ if and only if $\langle m,n\rangle=0\ \forall n\in N$
implies $m=0$, it is {\sl separated in} $N$ if and only if $\langle
m,n\rangle=0\ \forall m\in M$ implies $n=0$ and it is {\sl separated} if
and only if it is separated both in $M$ and in $N$. We already know that
if $\langle ,\rangle$ is separated in $M$, then $M$ is diagonal and if
$\langle, \rangle$ is separated in $N$ then $N$ is diagonal.\\ Finally a
central bimodule $M$ will be said to be {\sl reflexive} whenever
$M=M^{\ast_{A}\ast_{A}}$, which implies that $M$ is diagonal, and a
$Z(A)$-module $N$ will be said to be $A$-{\sl reflexive}, or simply {\sl
reflexive}, whenever $N=N^{\ast_{A}\ast_{A}}$, which implies that $N$ is
diagonal. If $M$ is reflexive then $M^{\ast_{A}}$ is reflexive and if
$N$ is reflexive then $N^{\ast_{A}}$ is reflexive.

\subsection*{Remark}

In fact the duality between central bimodules and $Z(A)$-modules comes
from a duality between bimodules and $Z(A)$-modules. Indeed, if $M$ is
an arbitrary bimodule over $A$, then $M^{\ast_{A}}=\hom^A_A(M,A)$ is
again canonically a module over the center $Z(A)$ of $A$. Furthermore
$M^{\ast_{A}\ast_{A}}=\hom_{Z(A)}(M^{\ast_{A}},A)$ is still a diagonal
bimodule and one has again a canonical bimodule homomorphism
$c_M:M\rightarrow M^{\ast_{A}\ast_{A}}$ which is, as a homomorphism of
$M$ onto $c_M(M)$, the functor Diag defined and studied in [8] and [9]
of the category of bimodules into the category of diagonal bimodules.
The very reason why we here restrict attention to central bimodules is
that only central bimodules reduce canonically to modules whenever $A$
is commutative. From the point of view of the above duality, the
diagonal bimodules and the $A$-diagonal $Z(A)$-modules are favoured and
of course, even more favoured are the reflexive bimodules and the
$A$-reflexive $Z(A)$-modules.\\ After having introduced a notion of
connection for central bimodules, it is natural to define a dual notion
for $Z(A)$-modules. Let $N$ be a $Z(A)$-module, {\sl a connection on}
$N$ {\sl related to} $A$, or simply {\sl a connection on} $N$ when no
confusion arises, is a linear mapping $\nabla, X\mapsto \nabla_X$, of
$\gder(A)$ into the linear endomorphisms of $N$ such that one has \[
\left\{ \begin{array}{ll} \nabla_{zX}(n) & = z\nabla_X(n)\\ \nabla_X(zn)
& = X(z)n + z\nabla_X(n) \end{array} \right.  \] $\forall n\in N,\
\forall X\in \gder(A)$ and $\forall z\in Z(A)$.\\ One defines, as in \S
2, the curvature $R$ of $\nabla$ by
$R_{X,Y}=[\nabla_X,\nabla_Y]-\nabla_{[X,Y]}$ and $R$ is now an
antisymmetric $Z(A)$-bilinear mapping of $\gder(A) \times \gder(A)$ into
the $Z(A)$-module $\hom_{Z(A)}(N,N)$. The set of connections on $N$ is,
if not empty, an affine space modelled on
$$\hom_{Z(A)}(\gder(A),\hom_{Z(A)}(N,N)).$$ The above definition is
justified by the following lemma.\\

\begin{lemma} Let $M$ be a central bimodule with a connection $\nabla$.
Then, there is a unique connection, again denoted by $\nabla$, on the
$Z(A)$-module $M^{\ast_{A}}$ which satisfies \end{lemma}
$$X(\mu(m))=\nabla_X(\mu)(m)+\mu(\nabla_X(m)),\ \forall X\in \gder(A),\
\forall \mu\in M^{\ast_{A}}, \forall m \in M.$$ {\sl Dually, let $N$ be
a $Z(A)$-module with a connection $\nabla$. Then there is a unique
connection, again denoted by $\nabla$, on the central bimodule
$N^{\ast_{A}}$ which satisfies}
$$X(\nu(n))=\nabla_X(\nu)(n)+\nu(\nabla_X(n)),\ \forall X\in \gder(A),\
\forall \nu\in N^{\ast_{A}}, \forall n \in N.$$

\noindent \underbar{Proof}. Define $\nabla_X(\mu)$ for $X\in \gder (A)$
and $\mu\in M^{\ast_{A}}$ by
$\nabla_X(\mu)(m)=X(\mu(m))-\mu(\nabla_X(m))$, then it is easy to show
that $\nabla_X(\mu)\in M^{\ast_{A}}$ and that $\nabla$ is a connection
on $M^{\ast_{A}}$ in the above sense. On the other hand $\nabla$ is
obviously unique under the condition of the lemma. The proof of the dual
statement is similar. $\square$\\ In the case where $M$ (resp. $N$) is
reflexive then the affine space of all connections on $M$ (resp. $N$)
and the affine space of all connections on $M^{\ast_{A}}$
(resp.$N^{\ast_{A}}$) are isomorphic under the above mapping.\\ More
generally, let $\langle , \rangle$ be a duality between a central
bimodule $M$ and a $Z(A)$-module $N$, then a pair $(\nabla,\nabla')$ of
a connection $\nabla$ on $M$ and a connection $\nabla'$ on $N$ will be
said to be {\sl compatible with the duality} $\langle, \rangle$ if one
has $X(\langle m,n\rangle)=\langle\nabla_X(m),n\rangle + \langle
m,\nabla'_X(n)\rangle,\ \forall X\in \gder(A),\ \forall m\in M$ and
$\forall n\in N$. If the duality is separated in $M$ (resp. $N$) then
given $\nabla'$ (resp. $\nabla$), if $\nabla$ (resp.  $\nabla'$) exists
it is unique.

\section{Derivations and forms}

As an illustration of the notions introduced in the latter section, let
us investigate the duality, between $\Omega^1_{\der}(A)$ and $\gder(A)$
and between $\gder(A)$ and $\underline{\Omega}^1_{\der}(A)$. We
summarize the result in the following theorem.

\begin{theorem} One has
$\underline{\Omega}^1_{\mer}(A)=(\Omega^1_{\mer}(A))^{\ast_{A}\ast_{A}}$.
More precisely, one has canonically:\\ a)
$\Omega^1_{\mer}(A))^{\ast_{A}}=\gmer(A)$ and the duality is
separated,\\ b) $(\gmer(A))^{\ast_{A}}=\underline{\Omega}^1_{\mer}(A)$
and the duality is separated.  \end{theorem}

\noindent \underbar{Proof}. By the universal property of $d:A\rightarrow
\Omega^1_{\der}(A)$, [8], we know that we have canonically
$\hom^A_A(\Omega^1_{\der}(A),M)=\gder(A,M)$ for any diagonal bimodule
$M$; so the equality of a) follows by taking $M=A$. The corresponding
duality is separated since $\Omega^1_{\der}(A)$ is diagonal (in fact
this follows directly from the definitions). On the other hand, the
equality b) is just the definition of $\underline{\Omega}^1_{\der}(A)$
and the corresponding duality is separated because a) implies that the
$Z(A)$-module $\gder(A)$ is $A$-diagonal. (Actually this last statement
also follows directly from the fact that if $X\in \gder(A)$ is such that
$\omega(X)=0,\ \forall \omega\in \underline{\Omega}^1_{\der}(A)$, then
$da(X)=X(a)=0,\ \forall a\in A$, which means $X=0$). $\square$\\ This
theorem shows exactly in what sense the minimal bimodule of
derivation-based one-forms $\Omega^1_{\der}(A)$ is ``dense" in the
maximal one $\underline{\Omega}^1_{\der}(A)$. Applied to the case where
$A$ is the Heisenberg algebra, it implies that the algebra
$\hat\Omega_{\der}(A)$ introduced in [5] in connection with the
noncommutative symplectic structure of quantum mechanics is just
$\underline{\Omega}_{\der}(A)$ (and in fact all the cochains in this
case).\\ In Section 4, we have defined a linear connection to be a
connection on $\underline{\Omega}^1_{\der}(A)$. Part b) of the theorem
shows that there is a more restrictive notion of linear connection,
namely a connection relative to $A$ on the $Z(A)$-module $\gder(A)$
because by applying the second part of lemma 1, to such a connection
corresponds a unique connection on $\underline{\Omega}^1_{\der}(A)$ and
this mapping is affine and injective. In fact, given a connection
$\nabla$ on $\gder(A)$ the torsion of the corresponding linear
connection can be identified with the element $T$ of
$\hom_{Z(A)}(\Lambda^2_{Z(A)}\gder(A),\gder(A))$ defined by
$$T_{X,Y}=\nabla_X(Y)-\nabla_Y(X)-[X,Y],\ \forall X,Y\in \gder(A).$$
Part a) of the theorem combined with lemma 1 shows that there is an even
more restrictive notion of linear connection, namely a connection on
$\Omega^1_{\der}(A)$.

\section{Reality and hermitian structures}

In this section $A$ is a unital $\ast$-algebra over $\Bbb C$. An {\sl
involutive bimodule} or a $\ast$-{\sl bimodule} over $A$ is a bimodule
$M$ equipped with an antilinear involution $m\mapsto m^\ast$ such that
$(amb)^\ast=b^\ast m^\ast a^\ast,\ \forall m\in M$ and $\forall a,b\in
A$. Dually an involutive $Z(A)$-module is a $Z(A)$-module $N$ equipped
with an antilinear involution $n\mapsto n^\ast$ such that
$(zn)^\ast=z^\ast n^\ast,\ \forall n\in N$ and $\forall z\in Z(A)$.
Given an involutive bimodule $M$ then the $Z(A)$-module $\hom^A_A(M,A)$
is an involutive $Z(A)$-module with involution $\mu\mapsto \mu^\ast$
given by $\mu^\ast(m)=(\mu(m^\ast))^\ast,\ \forall \mu\in \hom^A_A(M,A)$
and $\forall m\in M$. Given an involutive $Z(A)$-module $N$ then the
diagonal bimodule $N^{\ast_{A}}=\hom_{Z(A)}(N,A)$ is an involutive
bimodule with involution $\nu \mapsto \nu^\ast$ given by
$\nu^\ast(n)=(\nu(n^\ast))^\ast,\ \forall \nu\in N^{\ast_{A}}$ and
$\forall n \in N$. Elements of such sets satisfying
$\lambda=\lambda^\ast$ are called {\sl hermitian} or {\sl real}. The
$Z(A)$-module $\gder(A)$ is an involutive $Z(A)$-module with involution
$X\mapsto X^\ast$ defined by $X^\ast(a)=(X(a^\ast))^\ast,\ \forall X \in
\gder(A)$ and $\forall a\in A$. $\underline{\Omega}^1_{\der}(A)$ and
$\Omega^1_{\der}(A)$ are therefore involutive bimodules. More generally
one extends the involution to $\underline{\Omega}_{\der}(A)$ and
$\Omega_{\der}(A)$ by setting
$\omega^\ast(X_1,\dots,X_k)=(\omega(X^\ast_1,\dots,X^\ast_k))^\ast$ for
$\omega\in\underline{\Omega}^k_{\der}(A)$, (or $\Omega^k_{\der}(A))$ and
$X_i\in \gder(A)$. With this involution $\underline{\Omega}_{\der}(A)$
is a differential graded $\ast$-algebra in the sense that one has
$d(\omega^\ast)=(d\omega)^\ast$ and
$(\alpha\beta)^\ast=(-1)^{k\ell}\beta^\ast\alpha^\ast$ for $\omega\in
\underline{\Omega}_{\der}(A)$ and
$\alpha\in\underline{\Omega}^k_{\der}(A),\ \beta\in
\underline{\Omega}^\ell_{\der}(A)$; the subspace $\Omega_{\der}(A)$ is a
differential graded $\ast$-subalgebra.

\subsection*{Remark}

It is more or less well known that from the point of view of quantum
theory as well as from the point of view of spectral theory the good
generalization of the notion of algebra of real functions {\sl is not}
the notion of real associative algebra but is the notion of the real
Jordan algebra of all hermitian elements of an involutive complex
algebra, i.e. $\ast$-algebra, which plays the role of the noncommutative
generalization of the algebra of complex functions. It follows that what
must generalize the module of sections of a real vector bundle for
instance, or more generally the notion of module over an algebra of real
functions is not the notion of right or left module or a notion of
bimodules over a real noncommutative algebra but the set of real (i.e.
hermitian) elements of a $\ast$-bimodule over a $\ast$-algebra which
plays the role of the sections of the complexified vector bundle. Thus
the natural category at hand is the category of involutive central
bimodules over a $\ast$-algebra, and even more, if one thinks of real
vector bundles for instance, the category of involutive diagonal
bimodules and for the finite case the category of involutive reflexive
bimodules over a $\ast$-algebra, (with some other conditions replacing
projectivity). Notice also that one can alternatively use the dual
notion of the real elements of an involutive $Z(A)$-module or of an
involutive diagonal or involutive reflexive $Z(A)$-module. In fact,
there is a more restrictive notion of involutive diagonal and involutive
reflexive which we call {\sl diagonal involutive} and {\sl reflexive
involutive} which we now define. For any $\ast$-algebra $A$ and any set
$I$, $A^I$ is canonically an involutive bimodule. A {\sl diagonal
involutive} bimodule over $A$, (resp. a $A$-{\sl diagonal involutive}
$Z(A)$-module), is a $A$-bimodule (resp. a $Z(A)$-module) which is
isomorphic to an involutive subbimodule (resp. sub-$Z(A)$-module) of
$A^I$ for some set $I$. These notions are $A$-dual and therefore if $M$
is diagonal involutive $M^{\ast_{A}\ast_{A}}$ is also so, and if
furthermore $M=M^{\ast_{A}\ast_{A}}$ we say that $M$ is {\sl reflexive
involutive}. Notice that $\Omega_{\der}(A),\gder(A)$ and $\os_{\der}(A)$
are diagonal involutive.\\

\noindent Recall that a {\sl hermitian form} on a right $A$-module $E$,
[2], [3], is a sesquilinear mapping $h:E\times E\rightarrow A$ such that
$h(\varphi a,\psi b)=a^\ast h(\varphi, \psi)b$ and
$(h(\varphi,\psi))^\ast=h(\psi, \varphi),\ \forall \varphi,\psi\in E$
and $\forall a,b\in A$.\\ For a bimodule $M$, {\sl a right-hermitian
form on} $M$, or simply a {\sl hermitian} form on $M$ when no confusion
arises, will be a sesquilinear mapping $h: M\times M \rightarrow A$ such
that $h(ma,nb)=a^\ast h(m,n)b$ and $(h(m,n))^\ast=h(n,m),\ \forall
m,n\in M$ and $\forall a,b\in A$, as above, {\sl and} $h(m,cn)=h(c^\ast
m,n),\ \forall m,n\in M$ and $\forall c\in A$. The reason why the latter
condition has been included is that it allows to compose hermitian forms
on right modules with (right-) hermitian forms on bimodules. Namely if
$E$ is a right module with a hermitian form $h_E$ and if $M$ is a
bimodule with a right-hermitian form $h_M$ then one defines a hermitian
form $h$ on the right module $E\displaystyle{\otimesinf_A}M$ by setting
$h(\varphi\otimes m, \psi\otimes n)=h_M(m,h_E(\varphi,
\psi)n)(=h_M(h_E(\psi,\varphi)m,n)),\ \forall \varphi, \psi \in E$ and
$\forall m,n \in M$. It is also clear that if $E$ is a bimodule and if
$h_E$ is a right-hermitian form then the above definition gives a
right-hermitian form $h$ on the bimodule
$E\displaystyle{\otimesinf_A}M$. Furthermore, this composition of
(right-) hermitian forms is associative in an obvious sense. Assume now
that the positive cone $A^+=\{\sum_ia^\ast_ia_i\vert a_i\in A\}$ of $A$
is strict i.e. that one has $A^+\bigcap(-A^+)=\{0\}$, then a (right-)
hermitian form $h$ on a right module or a bimodule $E$ is {\sl positive}
if $h(\varphi,\varphi)\in A^+,\ \forall \varphi\in E$ and {\sl strictly
positive} if furthermore $h(\varphi,\varphi)=0$ implies $\varphi=0$.\\
Let $M$ be an involutive bimodule and let $g$ be a bimodule homomorphism
of $\displaystyle{M\otimesinf_A M}$ into $A$, i.e. $g\in
\hom^A_A(\displaystyle{M\otimesinf_A M},A)$, such that
$(g(m,n))^\ast=g(n^\ast,m^\ast)$ then $(m,n)\mapsto h(m,n)=g(m^\ast,n)$
is a right-hermitian form on $M$. Conversely, if $h$ is a hermitian form
on $M$ then one defines a $g\in \hom^A_A(\displaystyle{M\otimesinf_A
M},A)$ by setting $g(m,n)=h(m^\ast,n)$ and one has
$(g(m,n))^\ast=g(n^\ast,m^\ast)$. Such a
$g\in\hom^A_A(\displaystyle{M\otimesinf_A M},A)$ satisfying
$(g(m,n))^\ast=g(n^\ast,m^\ast)$ will be called a {\sl real inner
product} on the involutive bimodule $M$; $g(m,m)$ is real whenever $m$
is real. We shall say that $g$ is {\sl positive} (resp. {\sl strictly
positive}) whenever the corresponding hermitian form is so.\\ Let $M$ be
a bimodule and let $M'=\hom^A(M,A)$ be the left $A$-module dual of $M$
as a right $A$-module. The left module $M'$ is in fact a bimodule if one
defines $\alpha.a$ for $\alpha\in M'$ and $a\in A$ by
$(\alpha.a)(m)=\alpha(am)$, $\forall m\in M$. If $M$ is a central
bimodule, then $M'$ is also a central bimodule since, for $\alpha\in
M'$, $m\in M$ and $z\in Z(A)$, one has
$(z\alpha)(m)=z\alpha(m)=\alpha(m)z=\alpha(mz)=\alpha(zm)=(\alpha
z)(m)$. Assume now that $M$ is an involutive bimodule equipped with a
real inner product $g$. One defines a bimodule homomorphism $g^\sharp\in
\hom^A_A(M,M')$ by setting $g^\sharp(m)(n)=g(m,n)\ \forall m,n\in M$.
The real inner product $g$ on $M$ will be said to be {\sl nondegenerate}
whenever $g^\sharp$ is injective. If $g$ is strictly positive, then $g$
is nondegenerate.\\ Given an involutive central bimodule $M$, a
connection $\nabla$ on $M$ will be said to be {\sl real} if
$(\nabla_X(m))^\ast=\nabla_{X^{\ast}}(m^\ast)$. If $g$ is a real inner
product on $M$, a real connection $\nabla$ on $M$ will be said to be
{\sl compatible with} $g$ if one has
$$X(g(m,n))=g(\nabla_Xm,n)+g(m,\nabla_Xn),\ \forall m,n\in M, \forall
X\in \gder(A).$$ With obvious notations the above condition also reads
$$Xg(m\otimesinf_A n)=g(\nabla_X^{\otimes^{2}}(m\otimesinf_A n))\
\mbox{or}\ X\circ g=g \circ \nabla_X^{\otimes^{2}}.$$\\ Notice that a
nondegenerate real inner product $g$ on $\os^1_{\der}(A)$ is not yet a
complete noncommutative generalization of the notion of
pseudo-riemannian structure (and of riemannian structure whenever $g$ is
strictly positive); indeed the noncommutative generalization of the
symmetry is still missing.

\section{Noncommutative (pseudo-)riemannian\protect\newline structures}

In this section $A$ is again a unital $\ast$-algebra over $\Bbb C$. We
wish to investigate what kind of additional symmetry one has to impose
on a nondegenerate real inner product on $\os^1_{\der}(A)$ in order that
it can be considered as a noncommutative generalization of a
pseudo-riemannian metric. Although the solution is quite obvious in
simple situations, for instance if $A$ is finite-dimensional, this is
not the case for a general $\ast$-algebra $A$ as we shall see.
Fortunately, by taking a dual point of view, there is a natural
generalization of the notion of a pseudo-riemannian metric on the
$Z(A)$-module $\gder(A)$. We define a {\sl pseudo-metric} to be a
symmetric $Z(A)$-bilinear mapping $g_{\ast}$ of $\gder(A) \times
\gder(A)$ into $A$, i.e. $g_{\ast}\in (S^2_{Z(A)}\gder(A))^{\ast_{A}}$,
which is real, i.e. $(g_\ast(X,Y))^\ast=g_\ast(X^\ast,Y^\ast)$, and
which is nondegenerate in the sense that the corresponding mapping
$g_\ast^\sharp:\gder(A)\rightarrow \os^1_{\der}(A)$ defined by
$g^\sharp_\ast(X)(Y)=g_\ast(X,Y)$ is injective. A connection $\nabla$
relative to $A$ on the $Z(A)$-module $\gder(A)$ which is torsion-free,
i.e. which satisfies $\nabla_X(Y)-\nabla_Y(X)=[X,Y]$, and which is such
that one has
$Z(g_\ast(X,Y))=g_\ast(\nabla_Z(X),Y)+g_\ast(X,\nabla_X(Y))$ for
$X,Y,Z\in \gder(A)$ will be called a {\sl Levi-Civita connection} for
$g_\ast$. Summing over the cyclic permutations of the last equation with
signs $+ + -$ and using the symmetry and the vanishing of the torsion
one obtains $$ \begin{array}{ll} 2g_\ast(\nabla_X(Y),Z) = &
X(g_\ast(Y,Z))+Y(g_\ast(X,Z))-Z(g_\ast(X,Y))\\ &
+g_\ast([X,Y],Z)-g_\ast([Y,Z],X)+g_\ast([Z,X],Y).  \end{array} $$ So if
there exists such a Levi-Civita connection for $g_\ast$, then it is
unique since $g_\ast$ is nondegenerate. It follows from the reality of
$g_\ast$ and from the uniqueness that a Levi-Civita connection for
$g_\ast$ is real, i.e. that one has
$(\nabla_X(Y))^\ast=\nabla_{X^{\ast}}(Y^\ast)$. As pointed out in
Section 7, such a connection can be identified with a connection on
$\os^1_{\der}(A)$ (i.e. with a linear connection) which is torsion-free
and the above reality condition implies that it is a real connection on
$\os^1_\der(A)$ in the sense of Section 8. We are now in a position to
discuss the additional symmetry that one has to impose on a
nondegenerate real inner product on $\os^1_{\der}(A)$ in order that it
generalize a pseudo-riemmannian metric. Both
$\os^1_{\der}(A)\displaystyle{\otimesinf_A }\os^1_{\der}(A)$ and
$(S^2_{Z(A)}\gder(A))^{\ast_{A}}$ are sub-bimodules of the diagonal
bimodule $(\gder(A)\displaystyle{\otimesinf_{Z(A)}}\gder(A))^{\ast_{A}}$
of all $Z(A)$-bilinear mappings of $\gder(A)\times \gder(A)$ into $A$.
One defines a bimodule automorphism $\sigma$ of
$(\gder(A)\displaystyle{\otimesinf_{Z(A)}}\gder(A))^{\ast_{A}}$ by
setting $\sigma(b)(X,Y)=b(Y,X)$ for $b\in
(\gder(A)\displaystyle{\otimesinf_{Z(A)}}\gder(A))^{\ast_{A}}$ and
$X,Y\in \gder(A)$. The set of all $\sigma$-invariant elements
constitutes the bimodule $(S^2_{Z(A)}\gder(A))^{\ast_{A}}$ whereas
$\os^1_{\der}(A)\displaystyle{\otimesinf_A \os^1_{\der}}(A)$ is not
stable by $\sigma$ in general. The latter point is the only draw back to
writing the additional symmetry on the nondegenerate real inner product
on $\os^1_\der(A)$. Indeed, suppose that $A$ is such that
$\os^1_\der(A)\displaystyle{\otimesinf_A}\os^1_\der(A)$ is stable by
$\sigma$, for instance assume that
$$\os^1_{\der}(A)\displaystyle{\otimesinf_A}\os^1_\der(A)=(\gder(A)\displaystyle{\otimesinf_{Z(A)}}
\gder(A))^{\ast_{A}}$$ which is the case when $A$ is finite-dimensional,
then one can take the pseudo-metrics in
$\os^1_{\der}(A)\displaystyle{\otimesinf_A \os^1_{\der}}(A)$. One sees,
by duality, that in order that a nondegenerate real inner product $g$ be
a generalization of a pseudo-riemannian metric, it must be
$\sigma$-invariant, i.e. $g=g\circ \sigma$. In any case, in our
framework, we can content ourself with the above definition of
pseudo-metric. It is worth noticing that it has been suggested in [13]
that one can generalize our definition of linear connections in the case
where $\os_{\der}(A)\displaystyle{\otimesinf_A}\os_{\der}(A)$ is
$\sigma$-invariant to other differential calculi (non derivation-based)
by generalizing the bimodule homomorphism $\sigma$. This latter approach
has been used in two simple cases [7],[12].

\section*{Conclusion}

This paper is the first one of a series. Here we essentially introduce
the basic definitions and motivations without paying attention to the
existence problems. Also we have not introduced characteristic classes
but we have contented ourself with some comments on what they cannot be,
(factorization of inner derivations etc.). It must be clear that, in
order to define such classes as well as to develop a corresponding
$K$-theory, one must restrict attention to a class of bimodules (and
$Z(A)$-modules) which is smaller than the class of all central bimodules
(and all $Z(A)$-modules). It is also obvious that the (finite
projective) right and left modules together with their tensor products
and their tensor products with the appropriate bimodules have to be
taken into account. It is also worth noticing here that many notions
introduced in this paper do not refer to the specific differential
calculus (derivation-based) that we use and could be applied to other
differential calculi. Finally here we have worked in the purely
algebraic setting; but one can easily put everything in the setting of
convenient vector spaces in order to eventually take into account
topologies as in [9].  \section*{Acknowledgements}

It is a pleasure to thank John Madore and Thierry Masson for their kind
interest.

\end{document}